\documentclass[12pt,letterpaper]{amsart}
\usepackage{amssymb}
\usepackage{xy}
\xyoption{all}

\newtheorem{theo}{Theorem}[section]
\newtheorem{lemma}[theo]{Lemma}
\newtheorem{coro}[theo]{Corollary}
\newtheorem{prop}[theo]{Proposition}
\newtheorem{defi}[theo]{Definition}
\newtheorem{nota}[theo]{Notations}
\theoremstyle{definition}

\def\exp{\operatorname{exp}\nolimits}
\def\gcd{\operatorname{gcd}\nolimits}
\def\hd{\operatorname{hd}\nolimits}
\def\exp{\operatorname{exp}\nolimits}
\def\deg{\operatorname{deg}\nolimits}
\def\val{\operatorname{val}\nolimits}
\def\Disc{\operatorname{Disc}\nolimits}
\def\MaxSpec{\operatorname{MaxSpec}\nolimits}

\def\CL{{\mathcal{L}}}

\def\CH{{\mathcal{H}}}

\def\CR{{\mathcal{R}}}
\def\MJ{{\mathcal{J}}}

\def\CA{{\mathcal{A}}}

\def\CK{{\mathcal{K}}}
\def\CD{{\mathcal{D}}}
\def\tCD{\widetilde{\mathcal{D}}}

\def\BC{{\mathbf{C}}}
\def\BB{{\mathbf{B}}}
\def\bB{{\overline{E}}}
\def\Bb{{\mathbf{b}}}
\def\Bs{{\mathbf{s}}}
\def\BP{{\mathbf{P}}}
\def\BN{{\mathbf{N}}}
\def\BS{{\mathbf{S}}}
\def\BC{{\mathbf{C}}}

\def\BZ{{\mathbf{Z}}}

\def\tV{{\widetilde{V}}}
\def\tDelta{{\widetilde{\Delta}}}
\def\tf{{\widetilde{f}}}
\def\tW{{\widetilde{W}}}
\def\ts{{\widetilde{s}}}

\def\wt{\operatorname{wt}\nolimits}

\def\Id{\operatorname{Id}\nolimits}

\def\GL{\operatorname{GL}\nolimits}

\def\reg{\operatorname{reg}\nolimits}

\def\Feg{\operatorname{Feg}\nolimits}

\def\Spec{\operatorname{Spec}\nolimits}

\def\ie{{\em i.e.}}

\title{Zariski theorems and diagrams for
braid groups}
\author{David Bessis}
\thanks{The author thanks L\^e D\~ung Tr\'ang, Fabien Napolitano and
Bernard Teissier for useful conversations.}
\address{David Bessis, Department of Mathematics, Yale University,
P.O. Box 208283, New Haven CT 06520-8283, USA.}
\email{david.bessis@yale.edu}

\begin{document}
\begin{abstract}
Empirical properties of generating systems for complex reflection groups and
their braid groups have been observed by Orlik-Solomon
and Brou\'e-Malle-Rouquier, using Shephard-Todd classification. 
We give a general existence result for presentations of braid groups, 
which partially explains and generalizes the known empirical
properties.
Our approach is invariant-theoretic and does not use the classification.
The two ingredients are
Springer theory of regular elements and a Zariski-like theorem.
\end{abstract}
\maketitle    

{\bf Introduction}

Complex reflection groups share many properties with real reflection groups.
One of the main difficulties however is that no simple combinatorial
description of complex reflection groups (generalizing
Coxeter systems) is known. Elementary questions,
such as knowing how many reflections are needed to generate the group, do
not have satisfactory answers.
In \cite{orso}, Orlik and Solomon mention the following result
(where we have modified the notations to be consistent with the ones used here:
$d_1\leq \dots \leq d_r$ are the degrees, 
$d_1^* \geq \dots \geq d_r^*$ are the codegrees):

\begin{quotation}
{\small
{\bf (5.5) Theorem.} 
{\em Let W be a finite irreducible unitary reflection group. Then
the following conditions are equivalent:
\begin{itemize}
\item[($i$)] $d_i + d_i^* = d_r$ for $i= 1,\ldots,r$,
\item[($ii$)] $\sum_{i=1}^r (d_i+d_i^*)=rd_r$,
\item[($iii$)] $d_i^*<d_r$ for $i=1,\ldots,r$,
\item[($iv$)] $W$ may be generated by $r$ reflections,
\item[($v$)] If $\zeta=\exp(2i\pi / d_r)$ then there exist generating
reflections $s_1,\ldots,s_r$ for $W$ such that the element $c=s_1\ldots s_r$
has eigenvalues $\zeta^{d_1-1},\ldots,\zeta^{d_r-1}$ and
the element $c^{-1}$ has
eigenvalues $\zeta^{d_1^*+1},\ldots,\zeta^{d_r^*+1}$.
\end{itemize}
}
}
\end{quotation}
However, Orlik and Solomon describe these equivalences as ``surprising facts
for which we have no further explanation'', the proof relying entirely
on case-by-case study, using Shephard-Todd classification.
Our construction provides a partial ``explanation'' for these facts,
and reduces the use of the classification in the proof of the above theorem.

In \cite{brmarou} are listed diagrams which, for {\em almost all} irreducible
reflection groups,
symbolize presentations by generators and relations
of the associated braid group. These presentations satisfy remarkable
properties: the generators correspond to generators-of-the-monodromy,
the relations are positive and homogeneous
in the generators and, in many cases, a generator of the
center of the braid
group is given by a certain power of the product of all the generators,
taken in a certain order. 

The following theorem gives an {\em a priori}
invariant-theoretic explanation for the existence
of such presentations.
For simplicity, we 
state it only as an existence result, but as it will appear later
the proof of the existence of $\BS$ is constructive. The proof does not use the
classification.

\begin{theo}
\label{theointroduction}
Let $W$ be an irreducible complex reflection group of rank $r$,
with associated braid
group $\BB(W)$. Let $d$ be one of the degrees of $W$.
Assume that $d$ is a Springer regular number.
Let $N$ be the number of reflections in $W$, $N^*$ the number of
reflecting hyperplanes. Let $n:=(N+N^*)/d$. Then $n$ is an integer
and there exists a subset $\BS=\{\Bs_1,\dots,\Bs_n\}$ of $\BB(W)$, such that:
\begin{itemize}
\item[(a)] The elements $\Bs_1,\dots,\Bs_n$ are generators-of-the-monodromy,
and therefore their images $s_1,\dots,s_n$ in $W$ are reflections.
\item[(b)] The set $\BS$ generates $\BB(W)$,
and therefore $S:=\{s_1,\dots,s_n\}$
generates $W$.
\item[(c)] The product $(\Bs_1\dots\Bs_n)^d$ is central in $\BB(W)$ and
belongs to the pure braid group $\BP(W)$.
\item[(d)] Let $\zeta:=e^{2i\pi/d}$. 
The product $c=s_1\dots s_n$ is a $\zeta$-regular element in $W$ (and
therefore the eigenvalues of $c$ are $\zeta^{d_1-1},\dots,\zeta^{d_r-1}$, and
the eigenvalues of $c^{-1}$ are $\zeta^{d_1^*+1},\dots,\zeta^{d_r^*+1}$).
\item[(e)] There exists a set $\CR$ of relations of the form
$w_1=w_2$, where $w_1$ and $w_2$ are positive words of equal length in
the elements of $\BS$, such that $<\BS | \CR>$ is a presentation for 
$\BB(W)$.
\item[(f)] For $s\in S$, denote by $e_s$ the order of $s$. Take $\CR$ as
in {\em (e)} (but view it as a set of relations in $s_1,\dots,s_n$). The
group $W$ is presented by $<S | \CR \; ; \; \forall s\in S, s^{e_s}=1 >$.
\end{itemize}
\end{theo}

A consequence of this theorem is that {\em all} irreducible complex
reflection groups admit ``good diagrams''.
This is due to the fact (easily checked on the classification) that
any irreducible complex reflection group admits at least one regular degree.
Note that more than one degree may be regular. This may be an
indication that, taking into account the various complex symmetries of the
discriminant, one should perhaps consider more than one
``diagram'' for a single
reflection group. 

Sometimes the largest degree $d_r$ is regular; this is the case under
Orlik-Solomon assumptions $(i)$ or $(ii)$, and these assumptions imply
that the corresponding $n$ is $r$. Then $(iv)$ and $(v)$ are 
consequences of our theorem, which also implies the existence of ``good''
minimal presentations for the braid group, as in \cite{brmarou}.

For the exceptional $4$-dimensional group $G_{31}$, the largest degree
$24$ is still regular, but Orlik-Solomon assumptions $(i)$ or $(ii)$ are
not satisfied. However, we have analogs to $(iv)$ and $(v)$ where $r$
has to be replaced by the corresponding $n$ ($5$ in this example).
In addition, the theorem gives the existence (not previously known) of 
a Brou\'e-Malle-Rouquier-like presentation for the
braid group of $G_{31}$.

\section{Invariants, discriminant, regular degrees}

Let $V$ be $\BC$-vector space of finite dimension $r$.
A {\em (complex) reflection} of $V$ is a finite order element
$s\in\GL(V)$ such that $\ker(s-\Id)$ is an hyperplane.
A {\em reflection group} of $V$ is a finite subgroup $W\subset \GL(V)$
which is generated by reflections.

\subsection{Invariants} Let $W$ be a finite subgroup of $\GL(V)$.
Let $\BC[V]$ be the algebra of polynomial functions on $V$, \ie, the
symmetric algebra $S(V^*)$, with the natural grading (\ie, non
trivial linear forms are of degree $1$). We denote by $\BC[V]^W$ the subalgebra
of functions which are invariant for the dual action of $W$.
A classical theorem by Shephard-Todd (\cite{shephard}) states that
$\BC[V]^W \simeq \BC[V]$ if and only if $W$ is a reflection group.

From now on, we assume that $W$ is an irreducible reflection group.
Since $\BC[V]^W \simeq \BC[V]$, we may find $r$ algebraically independent
polynomial functions $f_1,\dots,f_r$ such that $\BC[V]^W=\BC[f_1,\dots,f_r]$.
We may actually require $f_1,\dots,f_r$ to be homogeneous, and
put in such an order that $\deg f_1 \leq \dots \leq \deg f_r$. Such a 
sequence $(f_1,\dots,f_r)$ is called a {\em system of basic invariants}.
The sequence $(d_1,\dots,d_r):=(\deg f_1,\dots,\deg f_r)$ is independent
of the choice of a specific system of basic invariants. The numbers
$d_1,\dots,d_r$ are the {\em degrees} of $W$. Note that they form
a multi-set, rather than a set, since there may be repetitions.

We denote by $\BC[X_1,\dots,X_r]$ the polynomial algebra in $r$ indeterminates.
A monomial $M=X_1^{\alpha_1}\dots X_n^{\alpha_n}$ has degree
$\deg(M)=\sum_{i=1}^n \alpha_i$. We also define its {\em weight}
by $\wt(M):=\sum_{i=1}^n d_i\alpha_i$.
This gives rise to two distinct graduations on $\BC[X_1,\dots,X_r]$. The one
associated with the degree is the {\em linear graduation}, the other
is the {\em weighted graduation}.

Choosing a system of basic invariant $f=(f_1,\dots,f_r)$ is the same
as choosing a graded algebra isomorphism $\Phi_f: \BC[X_1,\dots,X_r]
\stackrel{\sim}{\rightarrow} \BC[V]^W$, $X_i\mapsto f_i$, where 
$\BC[X_1,\dots,X_r]$ is endowed with
the weighted graduation.

Let $\BC[V]^W_+$ be the ideal of $\BC[V]^W$ consisting of elements vanishing
at $0$. A classical theorem of Chevalley (\cite{chevalley}) states that, as
a $\BC W$-module, $\BC[V]/\BC[V]^W_+\BC[V]$ is isomorphic to the regular
representation. It is a graded version of that representation, the grading
being inherited from the one on $\BC[V]$. 
 For $j\in \BZ_+$, denote by $(\BC[V]/\BC[V]^W_+\BC[V])_j$ 
the homogeneous summand of degree $j$.
For any irreductible $\BC W$-module
$M$, the {\em fake degree} of $M$ is the Poincar\'e polynomial of $M$
in $(\BC[V]/\BC[V]^W_+\BC[V])_j$:
$$\Feg_M(t) := \sum_{j=0}^{\infty} (M | (\BC[V]/\BC[V]^W_+\BC[V])_j ) t^j.$$
It is well known that
$$\Feg_V(t) = \sum_{i=1}^r t^{d_i-1}.$$
The numbers $d_i-1$ are called the {\em exponents} of $W$. By analogy,
Orlik and Solomon defined {\em coexponents}. We prefer here to work with
codegrees: the multiset of codegrees is the multiset $\{ d_1^*, \dots,
d_r^* \}$ uniquely defined by
$$\Feg_{V^*}(t) = \sum_{i=1}^r t^{d^*_i+1}$$
(the coexponents are the $d^*_i+1$). We will assume that the codegrees are 
ordered in {\bf decreasing} order:
$$d_1^* \geq \dots \geq d_r^*.$$

Let $\CA$ be the set of reflecting hyperplanes of $W$. For $H\in\CA$, we denote
by $e_H$ the order of the fixator $W_H$ of $H$ in $W$ ($W_H$ is the cyclic
group generated by the reflections with hyperplane $H$). Let $N$ be the number
of reflections in $W$, $N^*$ the number of hyperplanes. Obviously
$$ N = \sum_{H\in\CA} e_H-1 \quad \text{and} \quad N^* = \sum_{H\in\CA} 1.$$

The following lemma is classical.

\begin{lemma}
\label{NN*} We have
$$ N = \sum_{i=1}^r d_i -1 \quad \text{and} \quad N^* = \sum_{i=1}^r d^*_i +1.$$
\end{lemma}

\subsection{Discriminant}
The set of {\em regular vectors} in $V$ is 
$$V^{\reg} := V - \cup_{H\in \CA} H.$$
For each $H\in\CA$, choose a linear form $l_H\in V^*$ with kernel
$H$. An equation for $V^{\reg}$ in $V$ is $\prod_{H\in \CA}l_H \neq 0$.
The function $\prod_{H\in \CA} l_H^{e_H}$ is invariant, so
$\prod_{H\in \CA} l_H^{e_H}\neq 0$ is an equation for
$V^{\reg}/W$ in $V/W$.

\begin{defi}
The discriminant variety of $W$ is the closed subvariety of $V/W$ given
by the equation $\prod_{H\in \CA} l_H^{e_H} = 0$. 

The discriminant polynomial of $W$ with respect to a
system of basic invariants $f$ is
the multivariable polynomial
$\Delta_f := \Phi_f^{-1}(\prod_{H\in \CA} l_H^{e_H})$.
\end{defi}

It is easy to check that $\Delta_f$ is reduced (\ie,
$\BC[X_1,\dots,X_r]/ \Delta_f$ has no nilpotent elements).

The polynomial $\Delta_f$ is weighted homogeneous of weight $N+N^*$,
but is (in general) not
linearly homogeneous. Denote by $\deg(\Delta_f)$ the degree of a highest
degree monomial among those involved 
in $\Delta_f$, and by $\val(\Delta_f)$ the degree of the a
lowest degree monomial (such a monomial
is called a {\em valuation monomial}). The integer $\deg(\Delta_f)$ depends
on the choice of $f$, but 
$\val(\Delta_f)$ is independent from that
choice (the system-change morphisms $\Phi_{f'}^{-1}\Phi_f$ can only increase
the valuation and are invertible, therefore they preserve the valuation).
We have $$\val(\Delta_f) \geq \frac{N+N^*}{d_r},$$
since, for all monomial $M=X_1^{\alpha_1}\dots X_n^{\alpha_n}$,
we have $\wt(M)= \sum_{i=1}^n d_i\alpha_i \leq d_r \sum_{i=1}^r \alpha_i=
d_r \deg(M)$. 

Choose a base point $y_0$ in $V^{\reg}$, take its image $x_0$
as base point in $V^{\reg}/W$;
the pure braid group $\BP(W)$ and the braid group $\BB(W)$ associated with $W$
are, by definition, $\BP(W):= \pi_1(V^{\reg},y_0)$ and
$\BB(W):=\pi_1(V^{\reg}/W,x_0)$.
A (once again classical) theorem by Steinberg (\cite{steinberg})
states that
$V^{\reg}\twoheadrightarrow V^{\reg}/W$ is unramified, thus we have an 
exact sequence
$$\xymatrix@1{ 1 \ar[r] & \BP(W) \ar[r] & \BB(W) \ar[r] & W \ar[r] & 1.}$$

\subsection{Regular degrees}

The theory of regular elements has been initiated by Springer (\cite{springer}).

\begin{defi}
Let $\zeta\in\BC$ a root of unity. 
An element $w\in W$ is {\em $\zeta$-regular} if and only if 
$\ker(w-\zeta \Id) \cap V^{\reg} \neq \varnothing$.
An integer $d$ is a {\em regular number} for $W$ if and only if it is
the order of a regular element.
\end{defi}

Note that, by Steinberg theorem,
if $w$ is $\zeta$-regular, then $w$ and $\zeta$ have the same order.

\begin{lemma}
\label{prelang}
Let $f$ be a system of basic invariants. Let $d\in\BN$.
Let $\MJ$ be the ideal of $\BC[X_1,\dots,X_r]$ generated by
$\{X_j \; | \; d\nmid d_j \}$. Then
$$\text{$d$ is a regular number for $W$} \Leftrightarrow \Delta_f \notin \MJ.$$
\end{lemma}

\begin{proof}
Let $\zeta$ be a root of unity of order $d$.
Let $J$ be the subset of $\{1,\dots,n\}$ of those $j$ such that
$d \nmid d_j$.
By \cite{springer} 3.2 (i), the ideal in $\BC[V]$ of functions vanishing on
$\cup_{w\in W} \ker(w-\zeta \Id)$ is
$\sum_{j \in J} \BC[V]f_j$.
Saying that $d$ is regular is the same as saying that
$\prod_{H\in\CA} l_H^{e_H} \notin \sum_{j \in J} \BC[V]f_j$.
As $\BC[V]^W \cap \sum_{j \in J} \BC[V]f_j = \sum_{j\in J}\BC[V]^Wf_j$,
this is the same as saying that 
$\prod_{H\in\CA} l_H^{e_H} \notin  \sum_{j\in J}\BC[V]^Wf_j$, which, using
the isomorphism $\Phi_f$, is equivalent to 
$\Delta_f \notin \MJ$.
\end{proof}

Let $i\in\{ 1,\dots,r\} $. There is a canonical isomorphism
$\BC[X_1,\dots,X_r] \simeq \BC[X_1,\dots,\hat{X_i},\dots,X_r][X_i]$.
For any $P\in\BC[X_1,\dots,X_r]$, we denote by $P_{X_i}$
the one-variable polynomial with coefficients
in $ \BC[X_1,\dots,\hat{X_i},\dots,X_r]$ corresponding to $P$.

\begin{defi}
A polynomial $P\in\BC[X_1,\dots,X_r]$ is said to be
{\em monic in $X_i$} if and only if
the head coefficient of $P_{X_i}$ is a scalar.
\end{defi}

This notion makes it easy to recognize, among the degrees of $W$, the ones
which are regular:

\begin{lemma}
\label{lang}
Let $i_0\in \{ 1,\ldots,n\} $, let $d=d_{i_0}$.
\begin{itemize}
\item[(i)] The degree $d$ is regular if and only there exists a system of basic
invariants $f$ such that $\Delta_f$ is monic in $X_{i_0}$. When this
is the case, we have $\deg(\Delta_{f,X_{i_0}})= (N+N^*)/d$.
\item[(ii)] Assume that $\forall i \in \{ 1,\ldots,n\}, d | d_i \Rightarrow
i=i_0$. Then the following assertions are equivalent:
\begin{itemize}
\item  $d$ is regular,
\item there exists a system of basic
invariants $f$ such that $\Delta_f$ is monic in $X_{i_0}$,
\item for all system of basic
invariants $f$, $\Delta_f$ is monic in $X_{i_0}$.
\end{itemize}
\end{itemize}
\end{lemma}

\begin{proof}
(i): Let $f$ be a system of basic invariants.
Let, as in the proof of the previous lemma, $J \subset \{1,\ldots,n\} $ be
the subset of those $j$ such that
$d$ does not divide $d_j$, and $\MJ$ the ideal generated by
$\{X_j | j\in J \}$.
If $\Delta_f$ is monic in $X_{i_0}$, then $\Delta_f\notin \MJ$, so
$d$ is regular (by lemma \ref{prelang}).

Now assume that $d$ is regular, \ie, that $\Delta_f\notin \MJ$. Let
us construct from $f$ a system of basic invariants $f'$ such that
$\Delta_{f'}$ is monic in $X_{i_0}$.

Denote by $I$ the complement of $J$
(thus $i_0\in I$).
Let $a =(a_i)_{i\in\{ 1,\dots,n \}}$ be a family of complex numbers
such that $i\in J\cup\{ i_0\} \Rightarrow a_i=0$.
Let $f'_a:=(f_i - a_i f_{i_0}^{d_i/d})_{i\in \{ 1,\dots,n\}}$.
The requirements on $a$ ensure that $f'_a$ is a system of basic
invariants.
By replacing $X_i$ by $X_i + a_i X_{i_0}^{d_i/d}$ in $\Delta_f$, one obtains
the discriminant $\Delta_{f'_a}$.

Let $\overline{\Delta}_f$ be the image of $\Delta_f$ by the
composition $$\xymatrix@1{\BC[X_1,\dots,X_n] \ar@{>>}[r] & 
\BC[X_1,\ldots,X_n]/ \mathcal{J} \ar[rr]^{\sim}_{X_i+\MJ \mapsto X_i}
& & \BC[X_i,i\in I]}.$$ 
The polynomial $\overline{\Delta}_f$ is non-zero (by lemma \ref{prelang})
and weighted homogeneous of weight
$N+N^*$. There are two possibilities:
\begin{itemize}
\item Either $\overline{\Delta}_f \in \BC[X_{i_0}]$. By weighted homogeneity,
$\overline{\Delta}_f = b X_{i_0}^{(N+N^*)/d}$, with $b\neq 0$. Once again
by weighted homogeneity, $b X_{i_0}^{(N+N^*)/d}$ must be the only monomial of
 $\Delta_f$ of highest degree in $X_{i_0}$, and thus $\Delta_f$ is monic
in $X_{i_0}$.
\item Either $\overline{\Delta}_f \notin \BC[X_{i_0}]$. A direct
computation shows that the coefficient of $X_{i_0}^{(N+N^*)/d}$ in 
$\Delta_{f'_a}$ is the value of $\overline{\Delta}_f$ evaluated at
$X_{i_0}=1$ and $X_i=a_i, i\in I -\{i_0\} $. As 
$\overline{\Delta}_f \notin \BC[X_{i_0}]$, this coefficient, seen as a 
polynomial function with variables $a_i, i\in I-\{i_0\}$, is not constant. Thus
it is possible to choose the $a_i, i\in I-\{i_0\}$ such that this coefficient
is non-zero. By weighted 
homogeneity, the corresponding $\Delta_{f'_a}$ will be monic
in $X_{i_0}$.
\end{itemize}

(ii): The additional assumption, in the notations of the proof of (i),
is that $I=\{i_0\}$.
We prove that the first assertion implies the last one, which is enough.
Assume that $d$ is regular, let $f$ be a system of basic invariants.
The proof of (i) also proves that $\Delta_f$ is monic,
once it has been noticed that, in
the final discussion, $\overline{\Delta}_f \in \BC[X_{i_0}]$, since
$\BC[X_{i_0}]=\BC[X_i,i\in I]$.
\end{proof}

\begin{coro}
\label{postlang}
The largest degree $d_r$ is regular if and only if the valuation of
the discriminant is equal to $(N+N^*)/d_r$.
\end{coro}

\begin{proof}
If $d_r$ is regular, then by lemma \ref{lang} (i) we can find $f$
such that the monomial $X_r^{(N+N^*)/d_r}$ appears in $\Delta_f$; this
monomial must be a valuation monomial.

Now assume that the valuation of the discriminant is $(N+N^*)/d_r$. Choose
a system of basic invariants $f$. By weighted homogeneity, 
a valuation monomial in $\Delta_f$ can only
involve those $X_i$ for which $d_i=d_r$. Using
lemma \ref{prelang}, this implies that $d_r$ is regular.
\end{proof}

\section{A theorem \`a la Zariski}

Let $V$ an complex affine space, and
let $Z$ be an algebraic hypersurface of $V$.
Roughly speaking, Zariski theorems describe how the homotopy groups
of $V-Z$ can be compared to those of $H\cap(V-Z)$, where $H$
is a ``generic'' hyperplane. 

The situation we have in mind is when $Z$ is the discriminant
hypersurface of a complex reflection group, where we are looking for
generating systems for the fundamental group. A natural
way to map a free group to the fundamental group of $V-Z$ is by
considering an inclusion $L\cap(V-Z)\subset V-Z$, where $L$ is an
affine line. To obtain a surjective morphism, one could try to apply 
recursively a Zariski theorem, like the one in \cite{hammle}, to
a suitable generic affine flag. 

In our situation, certain directions are natural to consider:
the monic directions of the discriminant which, by lemma \ref{lang},
correspond to regular degrees. But there are no obvious ways to embed them in
complete flags. We state in this section a Zariski-like theorem
which allows us to work directly with an affine line, skipping the recursive
process; the genericity condition is explicit and elementary.
The proof follows Zariski's original strategy (see \cite{cheniot}
for a clear and modernized account).

\subsection{Generators-of-the-monodromy} First, we ought to
justify our use of dashes. To us, generators-of-the-monodromy are
just peculiar elements of fundamental groups; we do not want to discuss
what is the monodromy.
This subsection is included for the convenience of the reader; it is
almost ``copy-pasted'' from the appendix of \cite{brmarou}, where more
details can be found.

Let $U$ be a smooth connected complex algebraic variety. Let
$Z$ be an algebraic hypersurface of $U$, let $A$ be the set of irreducible
components of $Z$. Fix a basepoint $x_0\in U-Z$. A ``{\em path from
$x_0$ to $Z$ in $U-Z$}'' is a path $\gamma:[0,1] \rightarrow U$ such
that
\begin{itemize}
\item $\gamma(0)=x_0$,
\item $\gamma([0,1[)\subset U-Z$,
\item $\gamma(1)$ is a smooth point of $Z$. 
\end{itemize}
The reader should note (without being disturbed) that a path
from $x_0$ to $Z$ in $U-Z$ is not really a path in $U-Z$.
Let $\gamma$ be a path from $x_0$ to $Z$ in $U-Z$. A smooth point of $Z$
is a point which belongs to exactly one $D\in A$, and is smooth in that $D$.
Denote by $D_{\gamma}$ the ``target divisor'', \ie, the only $D\in A$
containing
$\gamma(1)$. We will say that ``{\em $\gamma$ is a path from $x_0$ to 
$D_{\gamma}$ (in $U-Z$)}'' (note that the notion really depends on the remaining
$D\in A$).

To $\gamma$, we associate the element $s_{\gamma} \in \pi_1(U-Z,x_0)$
which has the following informal description:
starting at $x_0$, follow $\gamma$; just before arriving at $\gamma(1)$, make
one full direct turn around $D_{\gamma}$; return to $x_0$ following
$\gamma$. The reader should check for himself that, as $\gamma(1)$ is 
a smooth point of $Z$, the ``local fundamental group'' of $U-Z$ at
$\gamma(1)$ is $\BZ$, and $s_{\gamma}$ is well-defined.
An element $s_{\gamma}$ obtained this way is a ``{\em generator-of-the-monodromy
around $Z$ in $U-Z$}'', or, more precisely, ``{\em around $D_\gamma$ in
$U-Z$}''.

Another path $\gamma'$ from $x_0$ to $D_{\gamma}$ is $D_{\gamma}$-homotopic
to $\gamma$ (relatively to $U-Z$)
if and only if there is a homotopy $\varphi:[0,1] \rightarrow
\Omega(U)$ such that $\varphi(0)=\gamma$, $\varphi(1)=\gamma'$ and
for all $t\in[0,1]$,
$\varphi(t)$ is a path from $x_0$ to $D_{\gamma}$. The element
$s_{\gamma}$ depends only on the $D_{\gamma}$-homotopy class of $\gamma$.

The following lemma is certainly well-known.

\begin{lemma}
\label{monodromy}
Let $U$ be a smooth connected 
complex irreducible variety, let
$A$ and $B$ be two families of irreducible codimension $1$
closed subvarieties. Assume that $A\cap B= \varnothing$.
Choose $x_0\in U -\cup_{D\in A \cup B} D$.
Consider the natural morphism
$\psi: \pi_1(U -\cup_{D\in A \cup B} D,x_0) \rightarrow
\pi_1(U - \cup_{D\in A} D,x_0)$. The morphism $\psi$ is surjective, and:
\begin{itemize}
\item[(i)] 
Let $D_0\in A$. For any generator-of-the-monodromy
$s$ around $D_0$ in $\pi_1(U - \cup_{D\in A} D,x_0)$, $\psi^{-1}(s)$
contains a generator-of-the-monodromy around $D_0$ in
$\pi_1(U -\cup_{D\in A \cup B} D,x_0)$.
\item[(ii)] 
The kernel of $\psi$ is the subgroup 
generated by the generators-of-the-monodromy
around $B$ in $\pi_1(U -\cup_{D\in A \cup B} D,x_0)$.
\end{itemize}
\end{lemma}

\begin{proof} As $U -\cup_{D\in A \cup B}D$ is obtained from
$U - \cup_{D\in A} D$ by removing complex codimension $1$ subvarieties,
$\psi$ is surjective.

In terms of paths, (i) is the following: let $\gamma$ be
a path in $U-\cup_{D\in A} D$ from $x_0$ to a divisor $D_{\gamma}\in A$;
then in the
$D_{\gamma}$-homotopy class of $\gamma$ relatively to $\cup_{D\in A} D$, there
exists a path avoiding $\cup_{D\in B} D$. This follows from standard general
position arguments, as $\cup_{D\in B} D \cap (U-\cup_{D\in A} D)$ has complex
codimension
$1$ in $U-\cup_{D\in A} D$ (whereas $\gamma$ has real dimension $1$),
and $\cup_{D\in B} D \cap D_{\gamma}$ has complex
codimension $1$ in $D_{\gamma}$ (whereas the target point $\gamma(1)$ has 
dimension $0$).

(ii) is nothing more than an induction from Proposition A.1 in \cite{brmarou}.
\end{proof}

Note that the assertion (ii) applied to $U=\BC^r$ and 
$A= \varnothing$ implies that the complement of an hypersurface
in $\BC^r$ is generated by (all) generators-of-the-monodromy.
In particular, braid groups are generated by generators-of-the-monodromy.

\subsection{The main tool}
Many Zariski-like theorems are either of projective (\cite{cheniot})
or local (\cite{hammle}) nature. The result proven here is truly of affine 
nature; it is both natural and
elementary, and we have been surprised not to find it in the litterature
(it may simply be that we did not look at the right place).

Let $\CH$ be an algebraic hypersurface in $\BC^r$, defined by a reduced
polynomial $P\in\BC[X_1,\dots,X_r]$. Choose $X=X_{i_0}$ one of the
indeterminates.
To simplify notations, we use $Y$ to refer collectively to the
variables $(X_i)_{i\neq i_0}$, e.g. we write $\BC[X,Y]$ instead
of $\BC[X_1,\dots,X_n]$.
A point in $\BC^r$ is described by its coordinates $(x,y)\in\BC\times\BC^{r-1}$.
Denote by $p$ the fibration $\BC^r\rightarrow\BC^{r-1}$, $(x,y)\mapsto y$.
The fibers of $p$ are the lines of
direction $X$; for $y\in\BC$, we denote by $L_y$ the line $p^{-1}(y)$.

\begin{nota}
We denote by the $\Disc(P_X)$ the resultant of $P'_X$ and $P_X$.
We denote by $\hd(P_X)$ the
head coefficient of $P_X$. We denote by $\alpha(P_X)$ the gcd of the
coefficients of $P_X$.
\end{nota}

Since $P$ is reduced, so is $P_X$, and $\Disc(P_X)\neq 0$.
We also have $$\alpha(P_X) = \gcd(P_X,\Disc(P_X)).$$

\begin{defi}
We say that a line $L_y$ of direction $X$ 
is {\em generic with respect to $\CH$} if and only
if $\Disc(P_X)(y)\neq 0$. A line $L_y$ of direction $X$ is {\em bad
with respect to $\CH$} if and only if $\alpha(P_X)(y)=0$.
\end{defi}

There are some lines which are neither generic nor bad.
The line $L_y$ is generic if and only if
the intersection $L_y\cap \CH$ has exactly $\deg (P_X)$ points. It is
bad if and only if $L_y\subset \CH$. The remaining lines are ``better'' than
generic, in the sense that $L_y\cap \CH$ is finite with cardinal strictly less
than $\deg(P_X)$.

Let $\CK$ be the hypersurface in $\BC^r$ defined by 
$\Disc(P_X)$. Let $\overline{\CK}$ be the hypersurface
in $\BC^{r-1}$ defined by $\Disc(P_X)$.
Let $$E:= \BC^r - (\CH\cup\CK), \quad \bB:= \BC^{r-1} - \overline{\CK}.$$
The restriction of $p$ makes $E$ into a fiber bundle over $\bB$, with fibers
being complex lines with $\deg(P_X)$ points removed.

Choose a basepoint $(x_0,y_0) \in E$, and take $y_0$ for basepoint in
$\bB$.

\begin{lemma}
\label{lift}
Let $s\in\pi_1(\bB,x_0)$ a generator-of-the-monodromy around $\overline{\CK}$.
Then there exists $\ts\in\pi_1(E,(x_0,y_0))$ a generator-of-the-monodromy
around $\CK$ such that $p_*(\ts)=s$.
\end{lemma}

\begin{proof}
The space $\BC^r-\CK$ is trivial bundle of fiber $\BC$ over
$\bB$, so $s$ can be lifted to a generator $s'$ of the monodromy around
$\CK$ in $\pi_1(\BC^r-\CK,(x_0,y_0))$ (any lifting of a defining path
suits). Let $A$ be the set of irreducible
components of $\CK$, let $B$ the set of irreducible components
of $\CH$ which are not in $\CK$. The point (i) of lemma \ref{monodromy}
applied to these $A$ and $B$ asserts that $s'$ is the image
of a generator-of-the-monodromy $\ts\in\pi_1(E,(x_0,y_0))$ around $\CK$ 
by the 
embedding morphism $\pi_1(E,(x_0,y_0)) \rightarrow \pi_1(\BC^r-\CK,(x_0,y_0))$. 
\end{proof}

Let $S\subset \pi_1(\BC^{r-1}-\overline{\CK},y_0)$
be a generating set of generators-of-the-monodromy around
$\overline{\CK}$ (this exists by lemma \ref{monodromy} (ii) ).
Using lemma \ref{lift}, 
lift $S$ to a set $\widetilde{S}\subset\pi_1(E,(x_0,y_0))$.
Using the fibration exact sequence 
$$\xymatrix@1{\cdots \ar[r]  &
\pi_1(L-L\cap\CH,x_0) \ar^{\phantom{mu}\iota_*}[r]
& \pi_1(E,(x_0,y_0)) \ar^{\phantom{mu}p_*}[r] &
\pi_1(\bB,y_0) \ar[r] & 1},$$
we see that  $\iota_*(\pi_1(L-L\cap\CH,x_0)) \cup \widetilde{S}$ generates
$\pi_1(E,(x_0,y_0))$.

Now make the additional assumption that $\alpha(P_X)=1$.
Let $A$ be the set of irreducible components of $\CH$, let
$B$ be the set of irreducible components of $\CK$.
The polynomials $P_X$ and $\Disc(P_X)$ are coprime, so
$A\cap B = \varnothing$.
By lemma \ref{monodromy} (ii), $\widetilde{S}$ belongs to the kernel
of the surjective morphism $\pi_1(E,(x_0,y_0)) \rightarrow 
\pi_1(\BC^r-\CH,(x_0,y_0))$.
This completes the proof of:

\begin{theo}
\label{theozariski}
Let $\CH$ be an algebraic hypersurface in $\BC^r$, defined by a reduced
polynomial $P\in\BC[X_1,\dots,X_r]$. Let $X=X_{i_0}$ be one of the
indeterminates.
Assume that the coefficients of $P_X$ are (all together)
coprime. Let $L$ be a line of direction $X$, generic with
respect to $\CH$. Then the inclusion
$L-L\cap \CH \hookrightarrow \BC^r - \CH$ is
$\pi_1$-surjective.
\end{theo}

The assumption
that the coefficients of $P_X$ are all together coprime is
satisfied, for example, when
$P_X$ is monic (this will be used
in the proof of theorem \ref{theointroduction}) 
or when $P_X$ is irreducible
(we will use this later, in section
\ref{drregular}).

\subsection{A refinement}
\label{refinement}
When $\alpha(P_X)=1$, the previous theorem yields generating sets
with $\deg(P_X)$ generators. 
Even when $\alpha(P_X)\neq 1$, 
``small'' generating sets can be constructed using the same strategy,
but the proof is slightly more technical. We begin by a definition:

\begin{defi}
Let $\sigma\in\mathfrak{S}_r$. 
We define an ordering $\leq_{\sigma}$
on the set of $r$-tuples
of integers: $(a_1,\dots,a_r)\leq_{\sigma} (b_1,\dots,b_r)$ if and only
if 
$$(a_{\sigma(1)},\dots,a_{\sigma(r)})\leq (b_{\sigma(1)},\dots,b_{\sigma(r)})$$
for the lexicographical order.

Let $P\in\BC[X_1,\dots,X_r]$.
A monomial $M=X_1^{a_1}\dots X_r^{a_r}$ involved in $P$ is
{\em $\sigma$-dominant}
in $P$ if and only if it is such that $(a_1,\dots,a_r)$ is maximal for
$\leq_{\sigma}$.
We say that $M$ is {\em dominant} in $P$ if and only if it is dominant for some
$\sigma\in\mathfrak{S}_r$.
\end{defi}

{\bf Example.} The dominant monomials have been underlined in the following
polynomial: $\underline{X_1^5X_2} + X_1^4X_2^2+ \underline{X_1X_2^3}
+ X_2^3$.

For the sake of simplicity, the next proposition is stated only as an existence
result, though the proof is actually constructive.

\begin{prop}
\label{refined}
Let $\CH$ be an
hypersurface in $\BC^r$ defined by
 a reduced polynomial $P$. Let $M=X_1^{a_1}\dots X_r^{a_r}$ be a monomial 
involved in $P$. Assume that $M$ is dominant.
Then the fundamental group of $\BC^r-\CH$ can be generated with $\deg(M)$
generators-of-the-monodromy.
\end{prop}

\begin{proof}
We describe an inductive construction procedure for such generating sets.
We keep the notations used in the previous subsection.
Let $\sigma\in\mathfrak{S}_r$ be such that $M$ is $\sigma$-dominant.

If $r=1$, then $M$ is the head monomial of $P=P_X$, $\alpha(P_X)=1$
 and theorem \ref{theozariski} gives generating sets with
 $\deg(P_X)=\deg(M)$ generators-of-the-monodromy.

Now assume $r>1$. Let $i_0:=\sigma(1)$, let $X:=X_{i_0}$, as in the
previous subsection. We have $\deg(P_X)=a_{\sigma(1)}$.
The monomial $M/X^{a_{\sigma(1)}}$ is dominant in
$\hd(P_X)$. Let $\CL$ be the hypersurface in $\BC^{r}$ defined
by $\hd(P_X)$, and let $\overline{\CL}:=p(\CL)$.
By the induction hypothesis, it is possible to generate
$\pi_1(\BC^{r-1}-\overline{\CL},y_0)$ by a set $T_0$ consisting of 
$\deg(M/X^{a_{\sigma(1)}})=\deg(M)-a_{\sigma(1)}$
generators-of-the-monodromy. As $\hd(P_X) | \Disc(P_X)$, we have an
exact sequence
$$\xymatrix@1{1\ar[r] & \ker \psi \ar[r]  & \pi_1(\BC^{r-1}-\overline{\CK},y_0)
\ar[r]^{\psi} & \pi_1(\BC^{r-1}-\overline{\CL},y_0) \ar[r] & 1}.$$
Using lemma \ref{monodromy} (i), lift $T_0$ to a same cardinality 
set $T$ of generators-of-the-monodromy around $\overline{\CL}$ 
in $\pi_1(\BC^{r-1}-\overline{\CK},(x_0,y_0))$. Let 
$U\subset \pi_1(\BC^{r-1}-\overline{\CK},y_0)$
be the set of all
generators-of-the-monodromy around the irreducible components of
$\overline{\CK}$ which
are not in $\overline{\CL}$. By lemma \ref{monodromy} (ii), we have
$\ker \psi = <U>$. Thus $\pi_1(\BC^{r-1}-\overline{\CK},y_0)$ is generated by 
$T\cup U$.

Using lemma \ref{lift}, lift $T$ and $U$ to same cardinality
sets $\widetilde{T}$ and $\widetilde{U}$ of generators-of-the-monodromy
around $\CK$ in $\pi_1(E,(x_0,y_0))$. Choose $L$ generic of direction
$X$. The fibration argument used before still proves that
$\iota^*(\pi_1(L-L\cap \CH,x_0)) \cup \widetilde{T} \cup \widetilde{U}$
generates $\pi_1(E,(x_0,y_0))$. As $\alpha(P_X)=\gcd(P_X,\Disc(P_X))$,
the common irreducible components of $\CH\cap\CK$ are irreducible
components of the hypersurface $(\alpha(P_X)=0)$,
thus, as $\alpha(P_X) | \hd(P_X)$, they are
also irreducible components of $\CL$. Using lemma \ref{monodromy} (ii),
this implies that elements of $\widetilde{U}$ are mapped to $1$ in 
$\pi_1(\BC^r-\CH,(x_0,y_0))$. To generate $\pi_1(\BC^r-\CH,(x_0,y_0))$,
it is enough to take the images in $\pi_1(\BC^r-\CH,(x_0,y_0))$
of a generating set of
$\pi_1(L-L\cap \CH,x_0)$ (of cardinal $\deg(P_X)=a_{\sigma(1)}$) and
of $\widetilde{T}$ (of cardinal $\deg(M)-a_{\sigma(1)}$).
\end{proof}

{\bf \flushleft Remarks.}
\begin{itemize}
\item The above refinement is still not optimal, since
$\alpha(P_X)$ may be reducible and may also 
strictly divide $\hd(P_X)$, thus some elements
of $\widetilde{T}$ may become trivial in the fundamental group of
$\BC^r-\CH$.
\item In the next sections, we will apply theorem \ref{theozariski}
and its refinement to situations where $P$ is weighted homogeneous.
In that case, $\pi_1(\BC^r-\CH)$ is isomorphic to the local 
fundamental group at $0$. Local Zariski theorems (\cite{hammle}) applied
to locally generic affine flags yield generating sets with
$\val(P)$ generators.
Discriminants of reflection groups happen to always
have valuation monomials which are dominant (this is not true of any weighted
homogeneous polynomial), so the two methods give the same number of
generators.
Our tool also works with 
directions which are not locally generic, e.g., for $X^2-Y^3$, we may
use either $X$ or $Y$.
\item The above constructions can of course be generalized using
algebraic (not necessarily linear) coordinate systems in $\BC^r$.
\end{itemize}

\section{Presenting braid groups}
\label{presenting}

We give in this section a constructive proof of the theorem stated in the
introduction. 
We return to the situation and notations of the first section, where $W$ is 
a reflection group in $\GL(V)$. 
Let $i\in\{1,\dots,r\}$. Let $d=d_i$, $X=X_i$. 
We make the following assumption:
\begin{quote}
\em The number $d$ is assumed to be regular for $W$.
\end{quote}

The construction of the generating set for $\BB(W)$ depends on
three successive choices:

{1- Choose a system of basic invariants such that $\Delta_f$ is 
monic in $X$}.
This is possible thanks to lemma \ref{lang} (i).
The space $V^{\reg}/W$ is isomorphic to the complement
of the hypersurface $\CH$ of $\BC^r$ defined by $\Delta_f=0$.

{ 2- Choose a generic line $L$ of direction $X$.}
By theorem \ref{theozariski}, the embedding
$L\cap (\BC^r - \CH) \hookrightarrow \BC^r - \CH$ is $\pi_1$-surjective.
The cardinality of the intersection of $L$ with $\CH$ is equal to 
the degree of $\Delta_{f,X}$ which, according to lemma \ref{lang} (i),
is equal to $n:=(N+N^*)/d$.

{ 3- Choose a basepoint $x_0\in L$ and a {\em planar spider}
$\Gamma$ from $x_0$ to $L\cap \CH$.}
What we mean by a {\em planar spider} from $x_0$ to $L\cap \CH$ is
a collection $\Gamma=\{\gamma_1,\dots,\gamma_n\}$
of non-intersecting (except of course at $x_0$)
paths in $L$ connecting $x_0$ to each of the points in $L\cap \CH$.
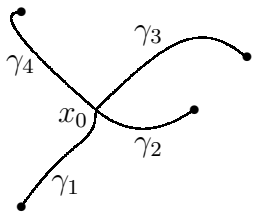
$$\xy
(-3,-1) *++{x_0}, (-10,13)="1", (0,0)="0", (20,7)="2",
(13,0)="3", (-10,-13)="4", "0";"1" **\crv{(-15,13)}, *\dir{*},
"0";"2" **\crv{(10,10)&(14,12)}, *\dir{*},
"0";"3" **\crv{(6,-5)}, *\dir{*},
"0";"4" **\crv{(0,-5)&(-3,-3)}, *\dir{*},
(-10,6) *++{\gamma_4},
(-4,-10) *++{\gamma_1}, (7,-5) *++{\gamma_2},
(7,10) *++{\gamma_3}
\endxy$$
(The suspicious reader may prefer to assume that all the paths
considered in this section are piecewise linear.)
We assume that the legs
$\gamma_1,\dots,\gamma_n$ are indexed in counterclockwise cyclic
order. 
To each leg $\gamma_i$, we associate the corresponding
generator-of-the-monodromy $\Bs_i \in \pi_1(\BC^r - \CH, x_0)$.
$$
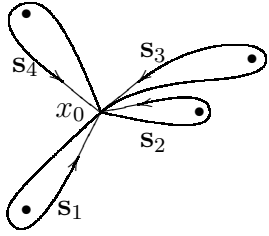
\xy
(-4,0) *++{x_0}, (-10,13)="1", (0,0)="0", (20,7)="2",
(13,0)="3", (-10,-13)="4", "0";"1" **@{}, *\dir{*},
"0";"2" **@{}, *\dir{*},
"0";"3" **@{}, *\dir{*},
"0";"4" **@{}, *\dir{*},
"0";(-5,4) **\crv{(-5,13)&(-10,17)&(-15,11)}, *\dir{<},
(-5,4);"0" **\crv{},
"0";(-3,-6) **\crv{(-12,-11)&(-14,-14)&(-9,-18)}, *\dir{<},
(-3,-6);"0" **\crv{},
"0";(5,1) **\crv{(13,-3)&(16,0)&(13,3)}, *\dir{<},
(5,1);"0" **\crv{},
"0";(5,4) **\crv{(5,4)&(20,4)&(23,7)&(20,10)&(10,8)}, *\dir{<},
(5,4);"0" **\crv{},
(-10,6) *++{\Bs_4},
(-4,-13) *++{\Bs_1}, (7,-4) *++{\Bs_2},
(7,8) *++{\Bs_3}
\endxy$$
Since $\BS:=\{\Bs_1,\dots,\Bs_n\}$ generates $\pi_1(L - (L\cup \CH),x_0)$,
it generates $\pi_1(\BC^r-\CH,x_0)$ which we identify with $\BB(W)$ through
the choice of $f$. Let $S:=\{s_1,\dots,s_n\}$ be the image of $\BS$ in
$W$. By \cite{brmarou} 2.14,
the image in $W$ of a generator-of-the-monodromy
is a reflection, thus the elements of $S$ are reflections. 

The properties (a) and (b) of \ref{theointroduction} are satisfied.
Let us now deal with (c) and (d) (note that the statement in (d) about
the eigenvalues is a standard property of regular elements, 
\cite{springer} 4.5).

\begin{lemma}
The product $(\Bs_1\dots\Bs_n)^d$ is central in $\BB(W)$, it belongs
to the pure braid group $\BP(W)$, and $c:=s_1\dots s_n$ is 
a $e^{2i\pi/d}$-regular element of $W$.
\end{lemma}

\begin{proof}
Note that $(\Bs_1\dots\Bs_n)^d$ central implies 
$$(\Bs_1\dots\Bs_n)^d = (\Bs_2\dots\Bs_n\Bs_1)^d = \cdots $$ and
that both being pure and being regular are invariant
by conjugation properties. Hence the
statement is independent from the choice of the starting leg $\gamma_1$
(only the cyclic order matters). This will allow us, later in the proof, 
to choose $\gamma_1$ at our convenience.

Let us check that it is enough to prove the statement for only one
basepoint chosen at our convenience.
Let $x_1$ be the basepoint in $L-L\cap\CH$ we prefer to the imposed $x_0$. To
compare the two situations, we can choose a path $\theta$ from
$x_0$ to $x_1$ in $L-L\cap\CH$. Along $\theta$ we have isomorphisms
$\pi_1(L-L\cap\CH,x_0)\simeq \pi_1(L-L\cap\CH,x_1)$ and
$\pi_1(\BC^r-\CH,x_0)\simeq \pi_1(\BC^r-\CH,x_1)$.
We can drag our original spider along $\theta$ to get a spider
with center $x_1$:
$$
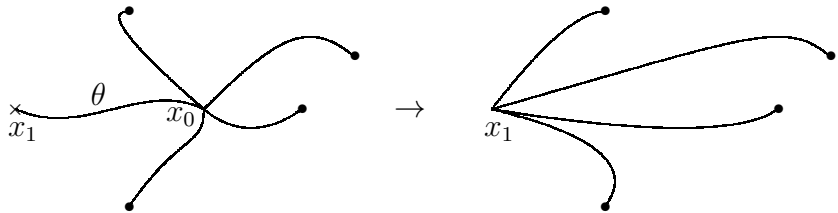
\xy
(-25,0)="x",
(-3,-1) *++{x_0}, (-10,13)="1", (0,0)="0", (20,7)="2",
(13,0)="3", (-10,-13)="4", "0";"1" **\crv{(-15,13)}, *\dir{*},
"0";"2" **\crv{(10,10)&(14,12)}, *\dir{*},
"0";"3" **\crv{(6,-5)}, *\dir{*},
"0";"4" **\crv{(0,-5)&(-3,-3)}, *\dir{*},
"0";"x" **\crv{(-8,4)&(-16,-4)}, *\dir{x},
(-14,2) *++{\theta}, (-24,-3) *++{x_1}
\endxy
\quad \rightarrow \quad
\xy
(-25,0)="x",
(-10,13)="1", (0,0)="0", (20,7)="2",
(13,0)="3", (-10,-13)="4", "x";"1" **\crv{(-15,13)}, *\dir{*},
"x";"2" **\crv{(10,10)&(14,12)}, *\dir{*},
"x";"3" **\crv{(6,-5)}, *\dir{*},
"x";"4" **\crv{(-4,-5)}, *\dir{*},
(-24,-3) *++{x_1}
\endxy
$$
The isomorphism
$\pi_1(\BC^r-\CH,x_0)\simeq \pi_1(\BC^r-\CH,x_1)$
maps the leg-generators of the original spider to the leg-generators
of the new one. It maps central elements, pure elements, and elements
whose image is regular to elements with the same property.
If the statement is proven for spiders with basepoint $x_1$, then
it follows for spiders with basepoint $x_0$.

Consider the affine segment $[0,x_0]$ in $\BC^r$.
For each $x\in [0,x_0]$, we denote by $L_x$ the affine line of direction $X$ 
passing through $x$ (we therefore have $L_{x_0}=L$).
As $\Delta_{f,X}$ is monic in $X$, $L_0\cap \CH = \{0\}$,
so $L_0$ is {\bf not} generic, unless $r=1$.
The lemma is obvious when $r=1$; we assume from now
on that $r>1$. This implies
that $L\neq L_0$ and that the segment $[0,x_0]$ has non-zero length and is
transverse to $X$. Let $E:=\cup_{x\in [0,x_0]} L_x$. The space $E$ is a
fiber bundle over $[0,x_0]$. The map $x \mapsto L_x\cap \CH$ from
$[0,x_0]$ to the space of
finite subsets of $E$ is continuous,
thus $E \cap \CH$ is compact. Let $R>0$ be large enough such that
$\forall x\in [0,x_0], L_x\cap\CH \subset B(x,R)$ (where $B(x,R)$ is the open
ball of center $x$ and radius $R$).

Take $x_1 \in L$ at distance $R$ from $x_0$. We will prove the statement
for spiders centered at $x_1$. 
$$
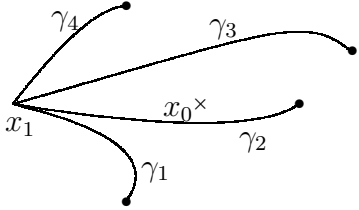
\xy
(-25,0)="x",
(-3,-1) *++{x_0},
(-10,13)="1", (0,0)="0", (20,7)="2",
(13,0)="3", (-10,-13)="4", "x";"1" **\crv{(-15,13)}, *\dir{*},
"x";"2" **\crv{(10,10)&(14,12)}, *\dir{*},
"x";"3" **\crv{(6,-5)}, *\dir{*},
"x";"4" **\crv{(-4,-5)}, *\dir{*},
"x";(0,0) **@{}, *\dir{x},
(-24,-3) *++{x_1},
(-18,11) *++{\gamma_4},
(-6,-9) *++{\gamma_1}, (7,-5) *++{\gamma_2},
(3,10) *++{\gamma_3}
\endxy
$$

To ease notations, let us fix affine
coordinates $$ \phi: E \stackrel{\sim}{\rightarrow} [0,1] \times \BC$$
such that $\phi(x_0)=(1,0)$ and $\phi(x_1)=(1,R)$. 
If the starting leg is well-chosen,
the product $\Bs_1\dots \Bs_n\in\pi_1(\BC^r-\CH,x_1)$
is represented by the loop
\begin{eqnarray*}
[0,1] & \longrightarrow & \BC^r-\CH \\
t & \longmapsto & \phi^{-1}((1,e^{2i\pi t} R))
\end{eqnarray*}

Let $x_2:=\phi^{-1}((0,R))$. The affine segment $[x_1,x_2]$ yields
an isomorphism $\pi_1(\BC^r-\CH,x_1) \simeq \pi_1(\BC^r-\CH,x_2)$.
As $R$ is large enough, the outer surface of the cylinder of radius $R$ around
$[x_2,x_1]$ does not intersect $\CH$, and the
isomorphism $\pi_1(\BC^r-\CH,x_1) \simeq \pi_1(\BC^r-\CH,x_2)$ maps
the product $\Bs_1\dots \Bs_n$ to the element
$\Bb$ represented by
\begin{eqnarray*}
[0,1] & \longrightarrow & \BC^r-\CH \\
t & \longmapsto & e^{2i\pi t} x_2
\end{eqnarray*}
(the $\BC$-action used in the formula is the linear one, not
the weighted action; in terms of weighted action, the formula
would of course be $t \mapsto  e^{2i\pi t/d} x_2$).
From that description of $\Bb$, it 
is classical (see for example \cite{brmi}, page 92) that $\Bb^d$ is central
in $\pi_1(\BC^r-\CH,x_2)$, that it is pure, and that
the image $b$ of $\Bb$ in $W$ is $e^{2i\pi/d}$-regular.
Thus $(\Bs_1\dots \Bs_n)^d$ is central and pure in $\pi_1(\BC^r-\CH,x_1)$,
and $c$, which is conjugate to the image $b$ of $\Bb$ in $W$, is a
$e^{2i\pi/d}$-regular element.
\end{proof}

The remainder of theorem \ref{theointroduction} can be deduced from
what we have obtained so far.
Let $F_{\BS}$ be the free group on $\BS$. Let
$\CR\subset F_{\BS} \times F_{\BS}$ be a set of relations
describing $\BB(W)$, \ie, such that the canonical morphism
$F_{\BS} \twoheadrightarrow \BB(W)$ is solution of the presentation
universal problem associated with $\CR$.
We have to prove that $\CR$ can be modified such that the relations
are in the free monoid on the alphabet $\BS$ and, in any relation,
the two sides have equal length.

Up to adding them (since, by the previous lemma, they are true), 
we may assume that $\CR$ contains the relations
$$(\Bs_1\dots \Bs_n)^d=(\Bs_2\dots \Bs_n\Bs_1)^d, \quad
(\Bs_2\dots \Bs_n\Bs_1)^d= (\Bs_3\dots \Bs_n\Bs_1\Bs_2)^d, \quad \cdots$$

There is a natural notion of length
in braid groups. 
Namely, the map $\BB(W) \rightarrow \BZ$ which maps
each generator-of-the-monodromy to $1$ extends to a unique morphism
$l: \BB(W) \rightarrow \BZ$. 
As each element of $\BS$ has length $1$, any relation $R\in \CR$
$$(R) \quad \quad \Bs_{i_1}^{\varepsilon_1}\Bs_{i_2}^{\varepsilon_2}\dots
\Bs_{i_k}^{\varepsilon_k} =
\Bs_{j_1}^{\varepsilon'_1}\Bs_{j_2}^{\varepsilon'_2}\dots
\Bs_{j_l}^{\varepsilon'_l},$$
(where the exponents are taken in $\{\pm 1\}$) must be homogeneous:
$$\sum_{p=1}^{k} \varepsilon_p = \sum_{q=1}^{l} \varepsilon'_q.$$
If any of the $\varepsilon_p$ or $\varepsilon'_q$ is negative, then by
multiplying both sides of $R$ by $(\Bs_1\dots \Bs_n)^d$, one gets an
equivalent relation in which the number of negative exponents has been
decreased (use the fact that $(\Bs_1\dots \Bs_n)^d$ is central thus
can be moved anywhere inside a word, and that it can be written
starting with any $\Bs\in\BS$). After a certain number of iterations, one
gets a relation $R'$ between positive words, which by homogeneity must have
equal length. The set $\CR':= (\CR-\{R\})\cup \{R'\}$ still describes $\BB(W)$.
By iteration, this proves property (e). Now (f) is a consequence of (e), using
Proposition 2.18 from \cite{brmarou}.

\section{Complements and applications}

\subsection{When the largest degree is regular}
\label{drregular}

The smallest generating sets obtained from
Theorem \ref{theointroduction} are obtained
with $d$ maximal among regular degrees. The ideal situation 
is when the maximal degree is regular.
By inspecting the classification, one sees that:

\begin{prop}
Let $W$ be an irreducible complex reflection group. Let $d_r$ be
the largest degree of $W$.
\begin{itemize}
\item If $W$ is a Coxeter group, or $G(d,1,r)$, or $G(2d,2,2)$, or
$G(e,e,r)$, or an exceptional group other than $G_{15}$, then
$d_r$ is regular. 
\item If $W$ is $G(de,e,r)$ with $e>2$, $d,r\geq 2$, or $G(2d,2,r)$ with
$d\geq 2$, $r>2$, or the exceptional group $G_{15}$, then $d_r$ is not
regular.
\end{itemize}
\end{prop}

The good news is that $d_r$ is regular for almost all exceptional groups,
including $G_{24}$, $G_{27}$, $G_{29}$, $G_{31}$, $G_{33}$ and $G_{34}$.
No presentations are known for the six corresponding braid groups. 
Theorem \ref{theointroduction} proves the existence of nice presentations.
Moreover, in each case, $n= (N+N^*)/d_r$ is the minimal number of
reflections needed, which is 
either $r$ (for $G_{24}$, $G_{27}$, $G_{29}$, $G_{33}$ and
$G_{34}$) or $r+1$ (for $G_{31}$), so the presentations are ``optimal''.

A criterion for deciding whether a number is regular or not has been
discovered by Lehrer and Springer (\cite{lehrerspringer}, 5.1):
$d$ is regular if and only if it divides as many degrees
as codegrees. Unfortunately, the ``if'' implication is proven by
case-by-case inspection (the ``only if'' is a consequence of elementary
properties of fake degrees).
Using theorem \ref{theointroduction} and Lehrer-Springer criterion,
one sees that, in the Orlik-Solomon theorem quoted in the introduction, 
assertion $(i)$ implies all the others.
Use of the classification is hidden in the Lehrer-Springer criterion.
A direct proof of this criterion, which is purely invariant-theoretic,
would be desirable.

As an example of a partial result which can be obtained without using
the classification (not even hidden in the Springer-Lehrer criterion),
let us mention the following. Let $W$ be an
irreducible complex reflection group of rank $r$. Assume that
the discriminant of $W$ is irreducible (this is equivalent, by classical
results, to the assumption that the $W$-action is transitive on the set
of reflecting hyperplanes; in other words, this is the analog for complex
reflection groups of the $ADE$ case for Weyl groups).
Then, in the Orlik-Solomon
theorem quoted in the introduction, assertion $(ii)$ implies all the others.
Indeed, since the discriminant $\Delta_f$ is irreducible, we can apply theorem
\ref{theozariski} (even without {\em a priori} knowing that $d_r$ is
regular) with a generic line of direction $X_r$. This yields a generating
set for the braid group of $W$ consisting of $\deg(\Delta_{f,X_r})$ 
generators-of-the-monodromy. Thus $W$ can be generated by
$\deg(\Delta_{f,X_r})$ reflections. Using the assumption $(ii)$ and 
weighted homogeneity, it is readily seen that $\deg(\Delta_{f,X_r}) \leq r$,
the equality being only possible when $\Delta_f$ is monic in $X_r$. 
Since an irreducible group of rank $r$ cannot be generated by less than
$r$ reflections, we see that $\deg(\Delta_{f,X_r}) = r$, which implies
assertion $(iv)$, and also that
$\Delta_f$ is monic in $X_r$. By lemma \ref{lang}, $d_r$ is regular.
We can now apply
theorem \ref{theointroduction}, which gives $(v)$. The assertions
$(i)$ and $(iii)$ are easily obtained by considering the eigenvalues of
a $d_r$-regular element of $W$ in the canonical representation and in the dual
representation.

\subsection{When the largest degree is not regular} As an example,
we discuss the case of the exceptional group $G_{15}$. The degrees
are $12$ and $24$, the codegrees $24$ and $0$. As $d_1=12$ is regular,
we may choose a $X_1$-monic discriminant polynomial:
$$X_1^5 + \alpha X_1^3 X_2 + \beta X_1X_2^2.$$
Up to replacing $X_2$ by $X_2 + \lambda X_1^2$, we may assume $\beta\neq 0$.
Applying theorem \ref{theointroduction} with $d=12$ yields a presentation
with $5$ generators. Using the refined Zariski theorem mentioned in
\ref{refined}, one gets a generating system associated with the
dominant monomial $\beta X_1X_2^2$, thus with $3$ generators, which is 
optimal.

The other cases can be handled the same way, as discriminants of 
complex reflection groups happen to always have a dominant valuation monomial.
More precisely, one can easily check on the classification
that if $W$ is an irreducible complex
reflection group of rank $r$
with non-regular highest degree, then, for a suitable
system of basic invariants $f$, the discriminant of $W$ has a
factorization $\Delta_f= X_i Q$, where $i\in\{1,\dots,n-1\}$ and $Q$ is monic
in $X_r$. 

\subsection{The minimum number of generators}
The situation can be summarized by the following proposition,
which partially follows from theorem \ref{theointroduction},
corollary \ref{postlang} and proposition \ref{refined}, and
partially from case-by-case analysis.

\begin{prop}
Let $W$ be an irreducible complex reflection group. The following integers
are equal:
\begin{itemize}
\item The minimum number of reflections needed
to generate $W$.
\item The minimum number of generators-of-the-monodromy
needed to generate $\BB(W)$,
\item The valuation of the discriminant.
\item $\lceil (N+N^*)/d_r \rceil$.
\end{itemize}
The degree $d_r$ is regular if and only if $(N+N^*)/d_r$ is an integer;
when this is the case, minimal presentations are described
by theorem \ref{theointroduction}.
\end{prop}

\subsection{Coxeter groups}

Let $W$ be an irreducible Coxeter group, seen as a complex
reflection group by complexifying the natural real reflection representation.
Examples of regular elements are
$w_0$ (the longest element) and $c$ (a Coxeter element). The
corresponding regular numbers are the degrees $2$ and $d_r$
($d_r$ is the Coxeter number, also denoted by $h$).
Generating sets corresponding to the standard Brieskorn presentation
for $\BB(W)$ can be obtained by applying theorem \ref{theointroduction}
with $d=d_r$.

Applying theorem \ref{theointroduction} with $d=2$ yields generating
sets with $N$ elements. An example of such a presentation is the
Birman-Ko-Lee presentation for $\BB(\mathfrak{S}_n)$ (\cite{BKL}).
Presentations of $\BB(W)$ with $N$ generators, with $W$ any irreducible
Coxeter group, are constructed in \cite{monoid}. They share with
Brieskorn and Birman-Ko-Lee presentations the following property:
the monoid presented with the same relations imbeds in the braid group.

Even when $W$ is not a Coxeter group,
any presentation obtained by theorem \ref{theointroduction} also
defines a monoid presentation, since the relations are positive.
The above embedding property is in general {\bf not} satisfied.

\subsection{Quotients and extensions of reflection groups}

In \cite{BBR} is studied a certain class of surjective morphisms between
reflection groups. Let us prove that, as announced in \cite{BBR},
the morphisms are induced from surjective morphisms between the corresponding
braid groups.

The situation is as follows. Let
$\tW$ be a reflection group acting on a complex vector space $\tV$.
Let $G$ be a normal subgroup of $\tW$. Let us assume that $G$ contains
no reflection. The action of $\tW$ on
$\tV$ induces an action on the variety $\tV/G$. Choose an imbedding
of $\tV/G$ in its tangent space $V$ at $0$. Let us assume that $G$ is
``good'' in $\tW$, as defined in \cite{BBR}, section 3.1.
This condition implies that the action on
$\tV/G$ is the restriction of an
action of $W:=\tW/G$ as a reflection group on $V$. Moreover, it is
possible to choose a system of basic invariants $f=(f_1,\dots,f_r)$ for $V$
(\ie, an identification $V/W \simeq \MaxSpec(\BC[X_1,\dots,X_r])$) and a subset
$I\subset \{1,\dots,r\}$ such that $\tf:=(\tf_i)_{i\in I}$ is
a system of basic invariants for $\tW$, where $\tf_i$ is the
composition $\tV\twoheadrightarrow \tV/G \hookrightarrow V 
\stackrel{f_i}{\rightarrow} \BC$.

Let $J:=\{1,\dots,r\} - I$.
The canonical embedding $$\Spec(\BC[X_1,\dots,X_r]/(X_j)_{j\in J})
\hookrightarrow \Spec(\BC[X_1,\dots,X_r])$$ identifies $\tV/\tW$ with
the linear subspace of $V/W$ defined by the equations $X_j=0$ for $j\in J$.
As observed in \cite{BBR}, paragraph 3.2.3, 
$\tV^{\reg}/\tW = (V^{\reg}/W) \cap \tV/\tW$.

Let us choose a basepoint $x_0\in\tV^{\reg}/\tW$. Define the braid
groups $\BB(\tW)$ and $\BB(W)$ with respect to $x_0$.

\begin{prop}
The natural morphism $\BB(\tW) \rightarrow \BB(W)$ is surjective.
\end{prop}

\begin{proof}
It is enough to prove the proposition when $\tW$ is irreductible:
if $\tW$ is reductible, it is a direct product of irreductible groups;
except in a few degenerate and straightforward cases, $G$ decomposes
as a corresponding direct product; then $W$ and $\BB(W)$ also 
decompose, and 
the reduction to the irreducible case follows.

Choose $d=d_{i_0}$ a regular degree for $\tW$ (one can check that
any irreducible reflection group admits at least one regular degree).

We may assume that the
discriminant $\tDelta_{\tilde{f}}$ is monic in $X_{i_0}$. Indeed, 
by lemma \ref{lang} (i), there
exists a system of basic invariants $\tf'$ of $\tW$ such that
the discriminant $\tDelta_{\tf'}$ of $\tW$ is monic in a variable
$X$. The system $\tf'$ is obtained from $\tf$ by a certain sequence of
algebraic substitutions. By performing the same substitutions among
the corresponding elements of $f$ (and leaving the remaining invariants
unchanged), we get a new system $f'$ which satisfies the 
monicity assumption (in
addition to the defining properties of $f$).

The discrimant $\tDelta_{\tf}$ of $\tW$ is obtained from the discriminant
$\Delta_f$ of $W$ by ``forgetting'' all monomials involving one
or more of the inderminates $X_j,j\in J$ (this operation is 
the composition of $\BC[X_1,\dots,X_r] \twoheadrightarrow \BC[X_1,\dots,X_r]/
(X_j)_{j\in J}$ with $\BC[X_1,\dots,X_r]/(X_j)_{j\in J} \simeq 
\BC[(X_i)_{i\in I}]$). As $\tDelta_{\tf}$ is monic in $X_{i_0}$, and
using weighted homogeneity, it is
readily seen that $\Delta_f$ was already monic in $X_{i_0}$,
and that  $\tDelta_{\tf}$ and $\Delta_f$ have the same degree in
$X_{i_0}$.

Let $\widetilde{\CH}$ be the hypersurface in $\tV/\tW$ defined by
$\tDelta_{\tf}=0$. Let $\CH$ be the hypersurface in $V/W$ defined
by $\Delta_f=0$.
Let $L$ be a $\widetilde{\CH}$-generic line of
direction $X_{i_0}$ in $\tW/\tV$. 
The cardinal of $L\cap \widetilde{\CH}=L\cap \CH$ is
$\deg(\tDelta_{\tf,X_{i_0}})=
\deg(\Delta_{f,X_{i_0}})$, thus $L$ is generic relatively to $\CH$ in 
$V/W$. By theorem \ref{theozariski}, the inclusion
$L-L\cap\CH \hookrightarrow V/W - {\CH}$ is $\pi_1$-surjective.
As it factors through
$\tV/\tW - \widetilde{\CH} \hookrightarrow V/W -\CH $, the latter map is
also $\pi_1$-surjective.
\end{proof}

Let $\tCD$ be a diagram for $\BB(\tW)$, symbolizing a presentation by generators
and relations as in section \ref{presenting}, \ie, corresponding to the choice
of a generic line of regular direction and of a ``planar spider''. The final
argument from the above proof makes it clear that, by adding some relations
to $\tCD$,
one gets a diagram $\CD$ for $\BB(W)$. The two diagrams are compatible, \ie,
the generators associated to $\tCD$ are sent to those of $\CD$. This explains
why the interpretation in \cite{BBR} of $\tW \twoheadrightarrow W$ as a 
``morphism of diagrams'' is actually valid at the level of braid groups.

\end{document}